\documentclass[a4paper,twoside,11pt]{amsart}   
 
\usepackage{setspace}% for double spacing  
\usepackage[inner=2.4cm,outer=4cm,top=2.4cm,bottom=2.4cm,marginparwidth=90pt]{geometry} %the geometry package allows different page sizes, margins etc.  
  
\usepackage{times}

\numberwithin{equation}{section}   
 
\usepackage{citesort}  
 
\usepackage{framed}  
\usepackage{comment}

\DeclareMathOperator{\divergenz}{div}  
\DeclareMathOperator{\dist}{dist}
\DeclareMathOperator{\graph}{graph}
\newcommand{\J}{\mathcal J}  
\renewcommand{\S}{\mathcal S}  
\newcommand{\dx}{\,dx}  
\newcommand{\ds}{\,ds}  
  
\newcommand{\dH}{\,dh}  
\newcommand{\bv}[1]{\int_\Omega \sqrt{1+\left|D#1\right|^2}}  
\newcommand{\BV}{{\mathrm{BV}}}  
\newcommand{\bigR}{{\mathbb R}}  
\newcommand{\eps}{\varepsilon}  
  
\newcommand{\wone}{\gamma_1}  
\newcommand{\wtwo}{\gamma_3}  
\newcommand{\wthree}{\gamma_2}
\newcommand{\wfour}{\gamma_4}   
\newcommand{\wi}{\gamma_i}

\newtheorem{theorem}{Theorem}[section]  
\newtheorem{lemma}[theorem]{Lemma}  
  
\newtheorem{corollary}[theorem]{Corollary}

\begin{document}  
  
\title{The capillarity problem for compressible liquids}  
\def\anuhome{@maths.anu.edu.au}
\def\monhome{@sci.monash.edu.au}

\author{Maria Athanassenas}  

\address{Maria Athanassenas,   School of Mathematical Sciences,  
Monash University,  
Vic 3800  
Australia}

\email{Maria.Athanassenas\monhome}
  
\author{Julie Clutterbuck}  

\address{Julie Clutterbuck,  Mathematical Sciences Institute,
Australian National University,
ACT 0200 Australia}  
  \email{Julie.Clutterbuck\anuhome}
  
\onehalfspacing  
  
 \subjclass[2000]{49Q20, 76N10, 76B45}
% 49Q20 Variational problems in a geometric measure-theoretic setting
% 76N10 Fluid mechanics, Compressible fluids and gas dynamics, general, Existence, uniqueness, and regularity theory
%76B45 Fluid mechanics, Incompressible inviscid fluids,  Capillarity (surface tension) 
\keywords{capillarity, functions of bounded variation, compressible liquids}
  
\begin{abstract}  
In this paper we study existence and regularity of solutions   
to the capillarity problem for compressible liquids in a tube. We introduce   
an appropriate   
space of functions of bounded variation, in which the energy functional   
recently introduced  by Robert Finn can be defined.  We prove   
existence of a locally Lipschitz minimiser in this class. 
\end{abstract}  

%\begin{titlepage}
\maketitle
%\end{titlepage}

\section{Introduction}  
 Extensive work has been published on  the behaviour of capillary (liquid-air or liquid-liquid)   
interfaces when the liquid is assumed to be incompressible. As an authoritative  
introduction we refer to \cite{finn-bible} by Finn. Two major approaches have been used to obtain existence and regularity results:  
classical PDE techniques for surfaces of prescribed mean curvature (see, for example, \cite{GT,huisken:1985,ladyzhenskaya-uraltseva}), and   
the functions of bounded variation and sets of finite   
perimeter setting for minimising the energy (see, for example, 
\cite{emmer,gerhardt:e-r,huisken:1985,giusti-pendent,tamanini-pendent}).
 
Results concerned with \emph{compressible} liquids are very recent and comparatively few,   
the model having been introduced by Finn in 2001
\cite{finn:2001}, see also \cite{finn-survey}.   
Following that paper we consider a capillary tube of cross section $\Omega\subset \bigR ^n$, which is simply connected and has Lipschitz boundary $\Sigma:=\partial\Omega$.  
We also assume that it satisfies an interior sphere condition of radius $R$.  
  
The capillary surface $\mathcal{S}$ is given as a graph of a function $u$ over the domain  
$\Omega$. We assume uniform downwards gravity $g$ and consider a {compressible} fluid of density $\rho$.  (In the incompressible fluid case, $\rho$ is constant.)

One can assume prescribed mass $M$,  
but the results in the present paper are for an infinite container.

We consider the energy for a capillary surface to consist of the following components:    

\emph{Energy of the free surface (surface tension):}
\begin{equation*}
E_{\mathcal{S}}=\frac\sigma{\rho_0}\int_\Omega \Phi(u;p_0)\sqrt{1+|Du|^2} \dx; 
\end{equation*}
  
\emph{Potential energy:}    
\begin{equation*}  
W=g\int_{\Omega}\int_0^u h\Phi(h;p_0) \dH \dx;  
\end{equation*}  
  
\emph{Wetting energy:}   
\begin{equation*}   
E_\Sigma=-\sigma\int_\Sigma \beta \int_0^u \Phi(h;p_0) \dH \ds;  
\end{equation*}  
here $\beta\in L^\infty(\Sigma)$  is the relative adhesion coefficient,  
satisfying $|\beta|\le 1-a$ with $a>0$; $\sigma$ and $g$ are the surface tension and gravitational constants; $\Phi(h;p_0)$ is the density function depending on
height $h$ and pressure $p$, which we assume to be given by one of the two models proposed by Finn \cite{athanassenas-finn,finn:2001}.
In the following, $p_0$ and $\rho_0$ will denote pressure and density at a reference level $u \equiv 0$.

\emph{Mass:}  
In the case of a mass constraint, a term $\lambda M$ is added to the energy, 
 where $\lambda$ is a Lagrange multiplier and the mass is   
\begin{equation*}  
M= \int_{\Omega}\int_0^u \Phi(h;p_0) \dH \dx.  
\end{equation*}  
  
The total energy (and in particular the wetting energy) need not be positive. 
  
A \emph{smooth} minimizer of the total energy $E_{\mathcal S}+W+E_\Sigma+\lambda M$ will satisfy the Euler-Lagrange equation
 \begin{equation}\label{cap surface equation}  
\divergenz {\frac{Du}{\sqrt{1+|Du|^2}}}= \frac{g{\rho_0}}\sigma  u +   \frac{D_1\Phi(u;p_0)}{\Phi(u;p_0)} \frac1{\sqrt{1+|Du|^2}} + \lambda \frac{\rho_0}\sigma   \,
\text{ on } \Omega,  
\end{equation}
with boundary condition
\begin{equation*}
\beta= \frac1{\rho_0}\frac{Du\cdot \nu}{\sqrt{1+|Du|^2}} \,\text{ on }\Sigma\label{cap surface boundary condition}, \notag   
\end{equation*} via standard calculus of variations techniques.

The present paper is based on one of the  models proposed by Finn for an \emph{isothermal} fluid: the density is assumed to be linear in the pressure, from which one obtains that $\Phi(h;p_0)=\rho_0e^{-\chi g h}$, for some positive constant $\chi$.

We may assume that $\chi=1$, $ g =1$, $\rho_0=1$, $\sigma=1$; other values of these constants correspond to different weightings on the components of the energy (that is, our energy becomes $\wone E_{\mathcal S}+ + \wthree W+ \wtwo E_\Sigma  + \wfour \lambda M$ for $\wi>0$), and a scaling of the domain $\Omega$.

Then   
\begin{equation*}  
\Phi(h;p_0)=e^{-h},  
\end{equation*}  
and the diverse components of the energy are as follows:  
\begin{gather*} E_{\mathcal{S}}=\int_\Omega e^{-u}\sqrt{1+|Du|^2} \dx,
\\ 
W=\int_{\Omega}\int_0^u h e^{-h} \dH \dx= \int_{\Omega}\left[1-e^{-u}(1+u)\right]\dx, 
\\ 
E_\Sigma=\int_\Sigma \beta \int_0^u e^{-h} \dH \ds= - \int_\Sigma \beta \left(1-e^{-u}\right) \ds,  
\\ 
M= \int_{\Omega}\int_0^u e^{-h} \dH \dx=\int_{\Omega} \left(1-e^{-u}\right) \dx.  
\end{gather*} 
As we are dealing with the case of an infinite container, we choose $\lambda=0$.  Without loss of generality 
(this will be shown when necessary, in Lemma \ref{B(k) is small}) 
 we may set  $\wi\equiv1$, and then seek to minimize the energy   
\begin{equation} \label{energy for u} \notag \J (u):= E_{\S}(u)+W(u)+E_{\Sigma}(u). \end{equation}

The following results have been recently obtained for the capillarity problem of a compressible fluid. 

For slightly compressible fluids Finn \cite{finn:2001} introduced the model we are using here. In the case of a tube closed at 
the bottom he found the necessary condition on the mass
 for existence of a solution is $M < \rho_0 |\Omega| / \chi g$.

For a circular tube, Finn and Luli \cite{finn-luli}  
show that for any boundary contact angle $\gamma$ with $0\leq \gamma < \pi$
there is at least one symmetric solution of the problem, and that the height of this solution will lie above any prescribed level if $M$ is sufficiently large. If $\gamma \leq \pi/2$, the solution is unique among symmetric solutions with that mass.

Finn and Athanassenas \cite{athanassenas-finn} follow the classical PDE approach.   
They include the situation where on the right hand side of the prescribed mean curvature equation \eqref{cap surface equation}   the term $\displaystyle{\frac{\rho_0 g}{\sigma}} u$ is replaced by    
$\displaystyle{\frac{\rho_0 - \chi p_0}{\sigma \chi}} (e^{\chi g u}-1)$, and study the non-constrained case.  
The results vary  
depending on   
the regularity of the boundary of the domain: they obtain height and gradient   estimates and existence of smooth solutions for smooth   
domains, but only \emph{variational 
solutions}  for domains with Lipschitz boundaries. 
As with the incompressible case, they observe that existence of solutions depends discontinuously on the opening angle of the corners of the domain.
In the case with the alternative right hand side of \eqref{cap surface equation}, they  show non-existence whenever the domain is small,  
that is,  when $\displaystyle{\frac{\rho_0 - \chi p_0}{\sigma \chi}} |\Omega| > - |\Sigma| \beta$ (here $\beta$ is taken to be constant).

In the present paper we use functions of bounded variation techniques.

In Section \ref{section2}, we introduce $\BV$, the space of functions of bounded variation.  After a transformation of $u$, the weighted surface area term is well defined in $\BV$.  Transforming the remaining components of the energy gives us a new energy, 
 $\J_1$.

In Section \ref{section4}, we prove height estimates. In Lemma \ref{ASlemma} and Lemma \ref{ASlemma-v2} we have two Stampacchia type results  needed in our case. 

In Section \ref{section5}, we show that the energy functional is bounded from below,   
and that a minimising sequence for the energy functional is uniformly bounded   
in the $\BV$-norm. Existence then follows via the standard compactness theorem and   
by the lower semicontinuity of the functional.  
  
Finally, in Section \ref{section8} we show that there exists a locally Lipschitz minimiser.

\subsection*{Acknowledgements}{This paper was begun during the first author's visit to the Max-Planck-Institut f\"ur Gravitationsphysik, Potsdam, and the second author's appointment at the Freie Universit\"at Berlin.  It was completed at the Centre for Mathematics and its Applications, Australian National University, Canberra.  We thank these institutions  for their support and hospitality. }

\section{The energy in the isothermal case}  \label{section2} 
As in \cite{giusti:bv2} the $\BV$-seminorms are:  
\begin{equation*}  
\int_\Omega \sqrt{1+|Du|^2} = \sup\left\lbrace\int_\Omega g_{n+1} + u 
\divergenz_n g \dx : 
 g_i\in C^1_0(\Omega)\,\, \forall i=1,\dots,n+1, 
\sum_{i=1}^{n+1}{g_i}^2\le 1\right\rbrace.  
\end{equation*}  
and   
\begin{equation*}  
\int_\Omega |Du| = \sup\left\lbrace\int_\Omega u    
\divergenz_n g \dx :  g_i\in C^1_0(\Omega)\,\, \forall i=1,\dots,n, \sum_{i=1}^{n}{g_i}^2\le 1\right\rbrace,  
\end{equation*}  
where $\divergenz _n$ is the divergence of the first $n$ components, $\divergenz_n g = \sum_{i=1}^n D_i g_i $.    
  
One then defines the spaces $\BV(\Omega):=\lbrace u\in L^1(\Omega): \int_\Omega \sqrt{1+|Du|^2} < \infty\rbrace$, and $\BV^+(\Omega):=\left\lbrace u\in \BV(\Omega): u\ge 0 \text{ almost everywhere in }\Omega\right\rbrace$.

In the case that $u\in C^1(\Omega)$, the surface energy term   
$E_{\mathcal{S}}=\int_\Omega e^{-u}\sqrt{1+|Du|^2} \dx$ may be simplified by  
writing $v=e^{-u}$.   Then we can rewrite it as $\int_\Omega \sqrt{v^2+|Dv|^2} \dx$, which bears close resemblance to the integral investigated  by Bemelmans and Dierkes  \cite{dierkes-bemelmans}, which was  $\int_\Omega \sqrt{v+|Dv|^2/4} \dx$; see also \cite{dierkes-huisken}.    
  
The focus of our investigation now shifts to $v$, rather than $u$ itself.

Define \begin{equation*}  
\int_\Omega \sqrt{ v^2+|Dv|^2}
:= \sup\left\lbrace\int_\Omega v\left(g_{n+1} + \divergenz _n  g \right) \dx :  g_i\in C^1_0(\Omega) \, \, 
\forall i=1,\dots,n+1, \sum_{i=1}^{n+1}{g_i}^2\le 1\right\rbrace.  
\end{equation*}  
  
\begin{lemma} \label{equiv if v smooth} If $v$ is smooth, $\int_\Omega \sqrt{ v^2+|Dv|^2}=\int_\Omega \sqrt{ v^2+|Dv|^2} \dx$.    
\end{lemma}    

\begin{proof} We consider the 
test function $g^\eps=\chi_\eps \left[v^2+|Dv|^2\right]^{-1/2}\left(-Dv,v\right)$, where $\chi_\eps$ is a sequence of $C^\infty_0(\Omega)$ functions with $\chi_\eps\le 1$, converging to $\chi_\Omega$, the characteristic function of $\Omega$, in $L^1$.    
Then    
\begin{align*}   
\int_\Omega \sqrt{ v^   
2+|Dv|^2} & \ge \int_\Omega v (g^\eps_{n+1} + \divergenz_n g^\eps) \dx \\   
&= \int_\Omega g^\eps_{n+1} v - Dv\cdot g^\eps  \dx \\  
&= \int_\Omega \chi_\eps \sqrt{ v^2+|Dv|^2} \dx\\  
&\rightarrow  \int_\Omega  \sqrt{ v^2+|Dv|^2} \dx   
\end{align*}  
as $\eps\rightarrow 0$.   
The other direction is similar. \end{proof}

We note the following fact:  
 \begin{lemma} \label{lemma 2.2}The quantity $\int_\Omega \sqrt{ v^2+|Dv|^2}$ is finite exactly when $v$ is in $\BV(\Omega)$.    
\end{lemma}

\begin{proof} Suppose that $v$ is in $\BV(\Omega)$.  Then   
\begin{align*}  
\int_\Omega \sqrt{ v^2+|Dv|^2} &= \sup\left\lbrace\int_\Omega v\left(g_{n+1} + \divergenz _n  g \right) \dx :     
 g_i\in C^1_0(\Omega), \,\sum_{i=1}^{n+1}{g_i}^2\le 1  
\right\rbrace  
\\&\le \sup\left\lbrace\int_\Omega \left(|v|+1\right)g_{n+1} + v \divergenz _n  g  \dx :    
 g_i\in C^1_0(\Omega),\, \sum_{i=1}^{n+1}{g_i}^2\le 1  
 \right\rbrace   
\\&\le  \sup\left\lbrace\int_\Omega g_{n+1} + v \divergenz _n  g  \dx :    
g_i\in C^1_0(\Omega),\, \sum_{i=1}^{n+1}{g_i}^2\le 1\right\rbrace  
\\&\phantom{====}+  \sup\left\lbrace\int_\Omega g_{n+1}|v|   
\dx :    
 g_{n+1}\in C^1_0(\Omega), \quad {g_{n+1}}^2\le 1 \right\rbrace   
\\&=\int_\Omega \sqrt{1+|Dv|^2} + \| v \|_{L^1(\Omega)}  
\\&<\infty.  
\end{align*}  
On the other hand, if $\int_\Omega \sqrt{ v^2+|Dv|^2}<\infty$, then $v\in L^1(\Omega)$, since if not, we can take $g_i=0$ for $i<n+1$ and $g_{n+1}=\chi_\eps$ (where $\chi_\eps$ is as in Lemma \ref{equiv if v smooth}) so that  
\begin{equation*}   
\int_\Omega \sqrt{ v^2+|Dv|^2}\ge \int_\Omega |v| \chi_\eps \dx \longrightarrow \infty   
\end{equation*}  
as $\eps\rightarrow 0$, contradicting our assumption.  Finally we can check  
\begin{align*}  
\int_\Omega   \sqrt{1+|Dv|^2}  
&= \sup\left\lbrace\int_\Omega g_{n+1} + v    
\divergenz_n g +vg_{n+1}-vg_{n+1} \dx :  g_i\in C^1_0(\Omega), \,\, \sum_{i=1}^{n+1}{g_i}^2\le 1\right\rbrace   
\\&\le  \sup\left\lbrace\int_\Omega  v    
\divergenz_n g +vg_{n+1} \dx :  g_i\in C^1_0(\Omega), \,\, 
\sum_{i=1}^{n+1}{g_i}^2\le 1\right\rbrace  
\\&\phantom{====} + \sup\left\lbrace\int_\Omega g_{n+1}(1-v)    
\dx :  
g_{n+1}\in C^1_0(\Omega), {g_{n+1}}^2\le 1\right\rbrace   
\\&\le   
\int_\Omega \sqrt{ v^2+|Dv|^2} + \|v\|_{L^1(\Omega)} +|\Omega| \\  
&<\infty.  
\end{align*}  
\end{proof}

 \begin{corollary}  \label{corollary 2.3}
If, in addition to the above conditions,  $v$ is in $\BV(\Omega)$, and $v_k$ is a mollification of $v$, then   
\begin{equation*} \int_\Omega \sqrt{ v_k^2+|Dv_k|^2} \rightarrow  
\int_\Omega \sqrt{ v^2+|Dv|^2}.  
\end{equation*}  
\end{corollary}  
This may be proved in the same manner as Lemma A1 of \cite{gerhardt:e-r}.

\label{section3}

Under the transformation $v=e^{-u}$, the wetting energy is $E_\Sigma=-\int_{\Sigma}\beta[1-e^{-u}]ds=-\int_{\Sigma}\beta[1-{v}]ds$, where we consider $\left.v\right|_\Sigma$ as a trace of $v$.   
 As in \cite[Theorem 2.10]{giusti:bv2}  
if  $\Sigma$ is Lipschitz, each function in $\BV(\Omega)$ has a trace in $L^1(\Sigma)$.        
Furthermore,  
if $\Sigma$ also satisfies an interior sphere condition with radius $R$, then the following  estimate holds (see \cite{gerhardt:e-r}, Remark 2):  
\begin{equation} \label{boundary estimate from CG}  
\int_{\Sigma}|v| \ds \le \int_\Omega |Dv| + c_{R} \int_\Omega |v| \dx,  
\end{equation}  
where $c_{R}$ depends on $n$, $R$, and $\Sigma$.  

The integrand of $W$, the potential energy term, becomes $\int_1^v \ln h \dH$ and 
so the complete energy,  in the isothermal case, is  
\begin{equation}\label{defn of J1}  
\mathcal{J}_1(v)=\int_\Omega \sqrt{v^2+\left|Dv\right|^2}+ \int_\Omega \int_1^v \ln h \dH \dx   
-\int_\Sigma \beta[1-v] \ds.     
\end{equation}

Here we are reminded  
 of the energy studied by Claus Gerhardt in \cite{gerhardt:e-r}, which was   
\begin{equation*}  
\bv{v} + \int_\Omega \int_0^v H(x,h) \dH\dx -\int_\Sigma \beta v \ds   
\end{equation*}  
for $\beta\in L^\infty(\Sigma)$  
and $H$ satisfying the conditions  
(a) $\dfrac{\partial H}{\partial h}>0$, and  
(b) $H(x,h_0)\ge (1+ c)$, $H(x,-h_0)\le -(1+c)$ for some $h_0\ge 0$ and a given $c$.    
  The conditions imply that for large values of $|v|$, the potential energy term (the one involving $H$) is strictly positive and increasing at least linearly in $|v|$.    The current case is an improvement on this situation:   a strictly positive potential energy, which increases like $v\ln v$ for large values of $v$.

\section{Height bounds} \label{section4} 
In this section, we assume that $v$ minimises $\J_1$ in $\BV(\Omega)$ and seek height bounds.  Note that a bound from above on $v$ would correspond to a bound from below on $u$, while a strictly positive bound from below on $v$ corresponds to a bound from above on $u$.   

At the end of this section we show an easier way to find one-sided estimates
 in the cases of 
$\beta$ being either positive (for which we show $v$ bounded from above) or negative (for which we show $v$ bounded from below).

\subsection{Height bounds from above on $v$} 
 To estimate a minimiser $v$ from above, we follow  an approach   
similar to \cite{huisken:1985,gerhardt:e-r} leading to a Stampacchia iteration \cite{STAMPACCHIA}. We use the following variant of the original Stampacchia lemma:   
    
\begin{lemma} \label{ASlemma}  
Suppose $B(t)$, non-negative and non-increasing in $t$, satisfies  
\begin{equation} (h-k)B(h)\le Ck \left[B(k)\right]^\gamma, \label{ASeqn}  
\end{equation}   
for  all $h,k$ such that $0<k_0\le k<h$, for some constants $C$, $k_0>0$ and $\gamma>1$.  Then $B(K)=0$ for some sufficiently large $K$ dependent on $C,\gamma,k_0$ and $B(k_0)$.    
\end{lemma}  
  
We will use the following in the proof of the above lemma:  
\begin{lemma} \label{convergence of sm} 
For all $\alpha>1$ and $d>-\alpha$, the sequence   
\begin{equation*} s_m=\left(1+\frac d {{\alpha^m}}\right)\left(1+\frac d {{\alpha^{m-1}}}\right)\dots\left(1+\frac d {\alpha}\right),\end{equation*}  
converges to a  non-zero limit.  
\end{lemma}  
  
\begin{proof}  We examine the sequence $\lbrace \ln s_n\rbrace$,  writing each term as the partial sum $\sum_{j=1}^m \ln\left(1+ \frac d {\alpha^j}\right)$, and using the ratio test for the convergence of series:  
\begin{align*} \lim_{j\rightarrow\infty} \left[ \ln\left(1+ \frac d {\alpha^{j+1}}\right)\right]\left[\ln\left(1+ \frac d {\alpha^j}\right)\right]^{-1}   
&= \lim_{j\rightarrow\infty}\left[ \dfrac{\partial}{\partial j} \ln\left(1+ \frac d {\alpha^{j+1}}\right)\right]\left[ \dfrac{\partial}{\partial j} \ln\left(1+ \frac d {\alpha^j}\right)\right]^{-1} \\  
&=\lim_{j\rightarrow\infty} \left[\frac{d\ln \alpha}{\alpha^{j+1}+d}\right]\left[\frac{d\ln \alpha}{\alpha^{j}+d}\right]^{-1} \\  
&=\frac 1\alpha.   
\end{align*}  
As this series converges to some limit $L$, $\lbrace s_m \rbrace$ converges to $e^{L}>0$.  
\end{proof}  
  
\begin{proof}[Proof of Lemma \ref{ASlemma}]  
We begin by defining the sequence  $k_m:=k_0 s_m $, where $s_m$ is as in the preceding lemma, with $\alpha=2$ and $d=C\left[B(k_0)\right]^{\gamma-1}2^{\gamma/(\gamma-1)}>0$ (we assume here that $B(k_0)\not=0$, otherwise the lemma is trivially true).  Note that as $k_{m+1} -k_m=k_0 d 2^{-(m+1)}>0$, $\lbrace k_m\rbrace$ is positive and increasing, and, by the above result, converges to some limit $K$.    
  
We now prove that $B(k_m)\le B(k_0)2^{\mu m}$ for $\mu=(1-\gamma)^{-1}<0$, by induction.    
  
The base step, for $m=1$, is as follows:  by assumption \eqref{ASeqn},  
\begin{equation*}  
(k_1-k_0)B(k_1)\le C k_0 B(k_0)^\gamma,  
\end{equation*}  
 and so using $k_1-k_0=k_0d /2$ we find that  
\begin{align*} B(k_1)&\le C   \frac{2^{1-\mu}}{d} B(k_0)^{\gamma-1} B(k_0)2^\mu\\  
&= 2^{1-\mu-\gamma/(\gamma-1)}  B(k_0)2^\mu\\   
&=  B(k_0)2^\mu.  
\end{align*}  
  
Now we make the inductive assumption that $B(k_m)\le B(k_0)2^{\mu m}$.   We use this and condition \eqref{ASeqn} to estimate  
\begin{align*}  
B(k_{m+1})&\le C \frac{2^{m+1}}{d} B(k_m)^{\gamma} \\  
&\le C \frac{2^{m+1}}{d} \left[B(k_0)2^{\mu m}\right] ^{\gamma}\\  
&= 2^{m+1-\gamma/(\gamma-1)}B(k_0)^{1-\gamma} \left[B(k_0)2^{\mu m}\right] ^{\gamma}  \\  
&\le B(k_0)2^{\mu(m+1)}.  
\end{align*}  
Finally, the monotonicity of $B$ implies that $B(K) \le \lim_{m\rightarrow \infty} B(k_m) \le \lim_{m\rightarrow \infty} B(k_0) 2^{\mu m}=0$.    
\end{proof}

\begin{theorem} \label{height bound from above on v}Let $v$ minimise $\J_1$ in $\BV^+ (\Omega)$, where  $\J_1$ is given by \eqref{defn of J1}.  Assume in addition $\partial\Omega$ to be Lipschitz and to satisfy an interior sphere condition. Then $v$ is bounded above.  
\end{theorem}  
\begin{proof}   

We set $A(k) = \{ x \in \Omega: v(x) > k \}$ for $k > k_0$, $k_0$ to be  
chosen later, the goal being to show that the non-increasing $|A(k)|$ vanishes for  
some large $k$. We also write $w := \min (v,k)$. As $v$ minimises $\mathcal{J}_1$,  
 we have $\mathcal{J}_1(v) \leq \mathcal{J}_1(w)$ for all eligible $w$, which   
after rearranging gives  
\begin{equation}  
0 \geq \left\{\int_{A(k)} \sqrt{v^2 + |Dv|^2}  - \int_{A(k)} k \dx  \right\}  
+ \int_{A(k)} \int_{k}^{v} \ln \, h \dH \dx + \int_{\Sigma} \beta (v-w) \ds.   
\label{v vs v-k equ}  
\end{equation}  
Here, we make use of the fact that $w \in \BV(\Omega)$, 
 that $Dw = Dv$ in $\Omega \setminus A(k)$ and $Dw = 0$ in $A(k)$.

We estimate the boundary term using \eqref{boundary estimate from CG}:  
\begin{equation*}   
 \left| \int_{\Sigma} \beta (v-w)  \ds \right| \leq (1-a) \left[ \int_\Omega  
\left| D(v-w)   \right| + c_{R} \int_{\Omega} \left| v-w  \right| \dx  
\right].  
\end{equation*}  
One can easily show that $\int_{\Omega} |Dv| \leq \int_{\Omega} \sqrt{v^2 + |Dv|^2}$,  
so that \eqref{v vs v-k equ} gives  
\begin{equation}  
k \left|A(k) \right| \geq a  \int_{A(k)} |D(v-k)| + \left(\ln \, k - (1-a)  
 c_{R}   \right) \int_{A(k)} \left|v-k \right| \dx.  
\label{Betweeneqn}  
\end{equation}  
  
For $\BV$ functions on $C^{0,1}$ domains $\Omega$, we have the Sobolev inequality  
$$\left[ \int_\Omega |f|^{\frac n{n-1}} \dx \right]^{\frac{n-1}n} \le c_2 \left[ \int_\Omega |Df| + \int_\Omega |f| \dx \right].$$ 
To see this, we first note that the inequality holds true for $W^{1,1}(\Omega)$ functions, since the space  
$W^{1,1}(\Omega)$ is continuously embedded in $L^{n/(n-1)}(\Omega)$ for $n>1$ (see Theorem 7.26 in \cite{GT}). The results extends to $f \in \BV (\Omega)$, after approximating $f$ by smooth functions  
as in Theorem 1.17 of \cite{giusti:bv2}, and then following the steps of the proof of Theorem 1.28 of  
\cite{giusti:bv2}.

We rearrange this inequality as  
$$ a  \int_\Omega |Df| \ge \frac a{c_\Omega} \left[ \int_\Omega |f|^{\frac n{n-1}} \dx \right]^{\frac{n-1}n} - a  \int_\Omega |f| \dx. $$    

Using the above estimate with $f=v-w$, \eqref{Betweeneqn} becomes   
\begin{equation*}   k|A(k)| \ge  
\frac a{c_\Omega} \left[ \int_\Omega |v-w|^{\frac n{n-1}} \dx \right]^{\frac{n-1}n} 
+ \left(\ln  k - (1-a)  
 c_{R}-a   \right)  \int_{A(k)} |v-w|   \dx.   
\end{equation*}   
By choosing $k_0$ (the lower bound on $k$) large enough, we can ensure that $\left(\ln  k - (1-a)  
 c_{R} -a  \right)$ is positive, and drop this term altogether.   The H\"older inequality \cite[Theorem 1.28]{giusti:bv2} gives    
\begin{equation*}
\| f\|_{L^{\frac n{n-1}}(A(k))}\ge \| f \|_{L^1(A(k))}|A(k)|^{-1/n},  
\end{equation*} and so for all $h,k$ with $h>k\ge k_0$, we have    
\begin{equation*}    
\frac a{c_\Omega}(h-k) |A(h)| \le \frac a{c_\Omega} \int_{A(k)} (v-k) \dx \le k|A(k)|^{1+\frac 1n}.   
\end{equation*}  
Now  
we can apply Lemma \ref{ASlemma} and conclude that for sufficiently large $K$, $|A(K)|=0$.   
\end{proof}  
  
\subsection{Height bounds on $v$ from below.}    
We start by remarking that bounds from below on $v$ --- or from above on $u$, where the surface is given by $\graph u$ --- are not essential for the existence proof. We will  see in the next section that the energy is bounded from below irrespective of such an estimate, and that subsequent results, leading to existence, also hold.  However, they are important for the correspondance between $u$ and $v$, and for the regularity results.

\begin{theorem} \label{bounds on u from above}  
Suppose that $v$ minimises $\J_1$ in $\BV^+(\Omega)$.    
Then there exists a bound from below on $v$, $0<c\le v$, $L^1$-almost everywhere.  
\end{theorem}  
The proof is similar to that of Theorem \ref{height bound from above on v}, but we use a slightly stronger Stampacchia-type result.  
\begin{lemma} \label{ASlemma-v2}  
Suppose $B(t)$, non-negative and non-increasing in $t$, satisfies  
\begin{equation} (h-k)B(h)\le C h\left[B(k)\right]^\gamma, \label{ASeqn-v2}  
\end{equation}   
for  all $h,k$ such that $0<k_0\le k<h$, for some constants $C$, $k_0>0$ and $\gamma>1$.   
If  
\begin{equation} \label{C small} 
CB(k_0)^{\gamma -1}<1,\end{equation} then there exists a $K<\infty$   
such that $B(K)=0$.    
\end{lemma}

\begin{proof}   
From \eqref{C small}, we may choose $d>0$ and $\alpha>1$ such that $CB(k_0)^{\gamma -1} \alpha^{\frac\gamma{\gamma -1}}\le d < \alpha$.   
  
Next, define the sequence  $k_m:={k_0}/s_m=k_{m-1}(1-\frac d {\alpha ^m})^{-1} $, where $s_m$ is as in Lemma \ref{convergence of sm} (but note the change of sign on $d$); here  $\lbrace k_m\rbrace $ is positive, increasing, and by Lemma \ref{convergence of sm}, convergent to some $K$.  
  
We now prove that  $B(k_m)\le B(k_0)\alpha^{\frac{-m}{\gamma-1}}$ by induction.    
  
The base step, for $m=1$, is as follows: by the definition of the sequence, we have  
 $k_1={k_0}(1-d/{\alpha})^{-1}$, and so our assumption \eqref{ASeqn-v2} gives  
\begin{equation*}(k_1-k_0)B(k_1)\le  C k_1B(k_0)^\gamma,\end{equation*}  
which leads to   
\begin{equation*} B(k_1)\le  \frac {C\alpha} d B(k_0)^{\gamma }   
\le B(k_0)\alpha^{\frac{-1}{\gamma -1}}.  
\end{equation*}   
  
Now we make the inductive assumption that $B(k_m)\le B(k_0)\alpha^{\frac{-m}{\gamma-1}}$, and show that this then holds for $k_{m+1}$:   we use  \eqref{ASeqn-v2} to estimate  
 \begin{align*}  
B(k_{m+1})&\le C \frac{k_{m+1}}{k_{m+1}-k_m} B(k_m)^\gamma    
\\&  
\le C \frac{\alpha^{m+1}}d B(k_0)^{\gamma-1} \alpha^{\frac{-m\gamma}{\gamma -1}+\frac{m+1}{\gamma -1}} \left[B(k_0)\alpha^{\frac{-(m+1)}{\gamma -1}} \right]  
\\&  
\le B(k_0)\alpha^{\frac{-(m+1)}{\gamma -1}}.   
\end{align*}  
  
Finally, the monotonicity of $B$ implies that $B(K) \le \lim_{m\rightarrow \infty} B(k_m) \le \lim_{m\rightarrow \infty} B(k_0) \alpha ^{\frac{-m}{\gamma -1}}=0$.    
\end{proof}  
  
We will need to show that the measure of the set where $v$ is small is small enough to satisfy \eqref{C small}.     This is the only place in this paper where it is not immediately clear that  rescaling the constants $\gamma_i$ to $1$ does not result in a loss of generality.   
Consequently, we include the arbitrary weightings in $\J_1$ in the following step.   
  
\begin{lemma}   \label{B(k) is small}  Let $v$ minimise $\J_1=\wone E_{\mathcal S}+\wthree W+\wtwo E_\Sigma$ in $\BV^+(\Omega)$, and set $B(k):= \lbrace x\in \Omega :v(x)<1/k \rbrace$.    Then for all $\eta>0$ we can find a $k$ such that $|B(k)|\le \eta$.    
\end{lemma}  
\begin{proof}  
Define the comparison function $w:=\max\lbrace v, \frac1k\rbrace\in \BV^+(\Omega)$ for any $k\ge k_0$, $k_0$ to be chosen later.   Note that $0\le w- v\le \frac 1 k$.    
Since $v$ minimises $\J_1$, we have $\J_1(v)\le \J_1(w)$.     
We use $\int\sqrt{u^2+|Du|^2}\ge \int |u| \dx$ for $u\in\BV^+(\Omega)$  
to  estimate  
\begin{align*}  
0&\ge \J_1(v)-\J_1(w)  
 \\&=   
\wone \int_{B(k)}\sqrt{v^2+|Dv|^2} - \wone \int_{B(k)} \frac1k \dx   
-\wtwo\int_{\Sigma} \beta(w-v)\ds   
\\&\phantom{++++}  
+\wthree \int_{B(k)} v\ln v-v -\frac 1k \ln \frac 1k +\frac1k \dx   
\\&\ge   
\int_{B(2k)} -|\wone-\wthree|\left|v-\frac 1k\right| +\wthree\left (v\ln v - \frac 1 k \ln \frac 1k\right) \dx  
\\&\phantom{++++}  
 +\int_{B(k)\setminus B(2k)} -|\wone-\wthree|\left|v-\frac 1k\right| +\wthree\left (v\ln v - \frac 1 k \ln \frac 1k\right) \dx  
\\&\phantom{++++}  
-\wtwo (1-a) |\Sigma| \frac 1 k   
\\ \intertext{--- now choose $k_0$ large enough so that $x\ln x$ is decreasing for $0<x\le \frac 1{k_0}$ ---}  
&\ge   
\int_{B(2k)} -|\wone-\wthree|\frac 1k +\wthree\left (\frac1{2k}\ln \frac1{2k} - \frac 1 k \ln \frac 1k\right) \dx  
 +\int_{B(k)\setminus B(2k)} -|\wone-\wthree|\frac 1{2k}  \dx  
\\&\phantom{++++}  
-\wtwo (1-a) |\Sigma| \frac 1 k   
\\&\ge   
-|\wone-\wthree|\frac 2 k |\Omega|  -\wtwo (1-a) |\Sigma| \frac 1 k + \wthree \int_{B(2k)} \left(-\frac 1{2k}\right) \left(\ln \xi +1\right) \dx \\  
\intertext{for some $\xi\in\left(\frac 1{2k},\frac 1k\right)$}  
&\ge  
-|\wone-\wthree|\frac 2 k |\Omega|  -\wtwo (1-a) |\Sigma| \frac 1 k   
+ \wthree |B(2k)| \left(-\frac 1{2k}\right) \left(\ln \frac 1{2k} +1\right).   
\end{align*}  
  Rearranging, and choosing $k_0$ large enough that $\ln\frac1{2k_0}<-1$, we find that   
\begin{equation*}  
|B(2k)|\le \frac{ 4 |\wone-\wthree|  |\Omega|  +\wtwo (1-a) |\Sigma| }  
{ \wthree \left(\ln {2k} -1\right)}<\eta   
\end{equation*}  
for  sufficiently large $k> k_0$.  
\end{proof}  
  
\begin{proof}[Proof  of Theorem \ref{bounds on u from above}]   
Let $B(k)$ be defined as above.   Set $w:=\max(v,1/k)$, for some $k\ge k_0$.  Again, as $v$ minimizes $\mathcal{J}_1$, then $\mathcal{J}_1(v)\le\mathcal{J}_1(w)$.  Proceeding exactly as in the proof of Theorem \ref{height bound from above on v}, we obtain  
\begin{align*}  
0&\ge a \int_\Omega |D(v-w)|  -\frac 1 k |B(k)| - \left[(1-a)c_1+ \ln(1/k)\right] \int_\Omega|v-w|\dx \\  &\ge \frac a {c_\Omega} \| v-w \|_{L^{\frac n{n-1}}(\Omega)}-\frac1 k |B(k)| - \left[(1-a)c_1+\ln(1/k)+a \right] \int_\Omega|v-w|\dx, \end{align*}  
and if we choose $k_0$ large, so that $\ln k\ge c_1(1-a)+a$, then the final term above is positive.  We drop it and apply the H\"older inequality to the $L^{\frac n {n-1}}$ term, leaving us with $(1/k)|B(k)|\ge (a/c_\Omega)\| v-w \|_{L^1(B(k))}|B(k)|^{-1/n}$, and so for each $h>k\ge k_0$ we have   
\begin{equation*}  
(h-k)|B(h)|\le C h |B(k)|^{1+\frac 1 n}.  
\end{equation*}  
  
Lemma \ref{B(k) is small} implies that we can find $k_0$ large enough that $C B(k_0)^{\frac 1 n }< 1$. 
We can then  
apply the Stampacchia-type Lemma \ref{ASlemma-v2}  to conclude that $|B(K)|=0$ for large   
$K$, and so $v\ge \frac1 K$ almost everywhere.    
\end{proof}  
 
\subsection{Height estimates in the cases $\beta\le0$ and $\beta\ge0$}  
Height estimates are easier to obtain in case $\beta$ is either non-positive or non-negative. 
 
We begin by observing a height bound for $v$ in the surface energy term.  This closely follows Lemma 5 of  
\cite{dierkes-bemelmans}, and may be proved in the same way.  
\begin{lemma}   
\label{lemma like lemma 5}  
Let $v\in\BV^+(\Omega)$ and suppose that $A(k)=\lbrace x\in\Omega: v(x)>k\rbrace$ has positive measure.  Then $w=\min(v,k)\in\BV^+(\Omega)$ and for almost all $k$,   
\begin{equation*}  
\int_\Omega \sqrt{ w^2+|Dw|^2}<\int_\Omega \sqrt{ v^2+|Dv|^2}. 
\end{equation*}\end{lemma}  
 
\begin{theorem} Suppose that $v\in\BV^+(\Omega)$ minimises $\mathcal{J}_1$, and that $\beta\ge 0$.  Then $v$ is bounded from above.  \label{height-bound-redundant-theorem}  
\end{theorem}  
\begin{proof} Set $w=\min(v,k).$  
  Suppose that $A(k)$ is of positive measure for some $k\ge 1$.  We may choose $k$ so that Lemma \ref{lemma like lemma 5} gives us 
\begin{equation*}  
\int_\Omega \sqrt{ w^2+|Dw|^2}-\int_\Omega \sqrt{ v^2+|Dv|^2}<0.
\end{equation*}  
We note that  
\begin{align*}  
\int_\Omega \int_1^w \ln h  \dH \dx -  \int_\Omega \int_1^v \ln h \dH \dx 
& =  \int_\Omega \int_v^w \ln h \dH \dx  
\\&=  \int_{\Omega \cap \lbrace x: v(x)\ge k\rbrace} \int_v^k \ln h \dH \dx 
\\&\le 0.    
\end{align*}  
Finally,   
\begin{equation*}  
-\int_\Sigma \beta[1-w] \ds+\int_\Sigma \beta[1-v]\ds= \int_\Sigma \beta [w-v] \ds \le 0  
\end{equation*}  
if $\beta\ge 0$.  Together, these inequalities give $\mathcal{J}_1(w)-\mathcal{J}_1(v)<0$, contradicting that $v$ was a minimum.  It follows that $|A(k)|$ cannot be positive, and so $v\le k$.      
\end{proof}

\begin{lemma}Suppose that $v$ minimizes $\mathcal{J}_1$, and $\beta\le 0$.  Then $v\ge e^{-1}$.    
\end{lemma}  
\begin{proof}   
Set $w=\max(v,\eps)$, and write $B(\eps)=\lbrace x\in\Omega: v(x)<\eps\rbrace$.  Then    
\begin{align*}  
\mathcal{J}_1(v)-\mathcal{J}_1(w)&= \int_{B(\eps)} \sqrt{ v^2+|Dv|^2}-\int_{B(\eps)} \sqrt{ w^2+|Dw|^2}- \int_\Omega \int_v^{\max(v,\eps)} \ln h \dH \dx   
\\&\phantom{spacespacespace}-\int_{\Sigma} \beta [w-v] \ds  
\\&=  
\int_{B(\eps)} \sqrt{ v^2+|Dv|^2}-\int_{B(\eps)} \eps \dx - \int_{B(\eps)} \int_v^{\eps} \ln h \dH \dx -\int_{\Sigma} \beta [\max(v,\eps)-v] \ds   
\\&\ge   
\int_{B(\eps)} \sqrt{ v^2+|Dv|^2}-\int_{B(\eps)} (\eps-v+v)\dx  + (-\ln{\eps})\int_{B(\eps)} (\eps-v) \dx   
\\&=  \int_{B(\eps)} \sqrt{ v^2+|Dv|^2}-\int_{B(\eps)}v\dx    
+ (-\ln{\eps}-1)\int_{B(\eps)} (\eps-v) \dx   
\\&>0  
\end{align*}  
for all $\eps<e^{-1}$, if $\int_{B(\eps)}|\eps-v|\not=0$.  However, this would contradict our assumption that $v$ is minimal for $\mathcal{J}_1$, so we conclude that $|B(\eps)|=0$ for small enough $\eps$.      
\end{proof}

\section{Existence of a minimiser}  \label{section5} 
 
\begin{lemma}[Lower bounds for the energy]  If $v\in\BV^+(\Omega)$, then $\J_1(v)\ge C(n,R,a,|\Omega|)$, where $C$ is not necessarily positive.   
\end{lemma}  
 
\begin{proof}  
As before, we can incorporate the wetting energy into the surface tension term using \eqref{boundary estimate from CG}, so that
\begin{align*}  
\J_1(v)&\ge \int_\Omega \sqrt{v^2+\left|Dv\right|^2}+ \int_\Omega \int_1^v \ln h \dH \dx   
-(1-a)c_R\int_\Omega|1-v| \dx -(1-a)\int_\Omega |Dv|   \\  
&\ge a\int_\Omega |Dv| + \int_\Omega f(v) \dx,  
\end{align*}   
where $f(v):= v(\ln v-1)+1-c_R(1-a)|1-v|$ is bounded below by a constant dependent on $c_R$ and $a$.   The result follows.  
\end{proof}  

 We  define a  \emph{minimising sequence for $\J_1$} as a   
sequence $v_j \in \BV^+ (\Omega)$ with   
$$\lim_{j\rightarrow\infty} \J_1(v_j)=\inf_{w\in\BV^+(\Omega)} \J_1(w):= m .$$

\begin{lemma} \label{sequence bounded in bv}  
A minimising sequence for $\J_1$ is uniformly bounded in the $\BV$-norm.  
\end{lemma}  
 \begin{proof}   
We can assume that  
$  
\J_1(v_j) \leq m+1  
$  
for $j$ large enough.  As in the previous lemma, where we defined $f$, we then have  
\begin{equation*}  
m+1\ge \J_1(v_j) \ge a\int_\Omega |Dv_j| + \int_\Omega f(v_j) \dx,  
\end{equation*}   
so the uniform bound follows from the lower bound on  $f$:  
\begin{equation*}  
\int_\Omega |Dv_j| \le \frac 1a \left( m+1 -  |\Omega|\inf_{h\in \bigR^+} f(h)\right).   
\end{equation*}  
Also, since there exist positive constants $\alpha_1$,$\alpha_2$ such that $f(t)\ge \alpha_1 t-\alpha_2$, we have the uniform $L^1$ bound 
\begin{equation*} 
\Vert v_j \Vert_{L^1(\Omega)}\le \frac1{\alpha_1}\left(\int_\Omega \left[f(v_j)+ \alpha_2\right] \dx\right)\le \frac1{\alpha_1}(m+1+\alpha_2|\Omega|). 
\end{equation*} 
\end{proof}  
   
\begin{lemma}[Lower semicontinuity of $\J_1$] \label{lsc for J1}  
A sequence $v_k\in\BV^+(\Omega)$ with  $v_k \rightarrow v$ in $L^1(\Omega)$  satisfies  
$$\J_1(v)\le  \liminf_{k\to\infty} \J_1(v_k).$$  
\end{lemma}  
\begin{proof}   
We show the surface energy term is lower semicontinuous.  For any admissible $g$, we have   
\begin{align*}\int_\Omega v \left(g_{n+1}+ \divergenz_n g\right) \dx   
&= \lim_{k\rightarrow\infty}\int_\Omega v_k \left(g_{n+1}+ \divergenz_n g\right) \dx \\   
&= \liminf_{k\rightarrow\infty}\int_\Omega v_k \left(g_{n+1}+ \divergenz_n g\right) \dx \\  
&\le \liminf_{k\rightarrow\infty} \int_\Omega \sqrt{ v_k^2+|Dv_k|^2}.  
\end{align*}  
Lower semicontinuity follows by taking the supremum over all admissible $g$.    
  
Continuity of the remaining terms of $\mathcal{J}_1$ follows as in \cite[Appendix II]{gerhardt:e-r}.  
\end{proof}

Combining all of the above results we have: 
\begin{theorem}[Existence of a minimiser]  
There exists a function $v \in \BV^+(\Omega)$, such that  
\begin{equation*}
\mathcal{J}_1(v) = \inf_{w \in \BV^+(\Omega)} \J_1(w).\end{equation*}
\end{theorem}  
 
\begin{proof}  
Let $\lbrace v_j\rbrace$ be the minimising sequence of Lemma \ref{sequence bounded in bv}, with $\Vert v_j \Vert_{\BV(\Omega)} \le C$.  By the standard compactness theorem (for example \cite[Theorem 1.19]{giusti:bv2})  there exists a subsequence $v_{j'}\to v$ in $L^1(\Omega)$.  
 
Since the $\BV$-norm is lower semicontinuous, $v$ is also in $\BV^+(\Omega)$,  and $\J_1(v)\ge \inf_{v\in \BV^+(\Omega)} \J_1(w)=m$.   Lower semicontinuity of $\J_1$, as in Lemma \ref{lsc for J1}, gives $\J_1(v)\le \liminf \J_1(v_{j'})=m$, completing the proof.  
\end{proof}

 \section{Regularity} \label{section8} 
In this section we show that a minimiser $v \in \BV^+(\Omega)$ of $\J_1$
is locally Lipschitz in
$\Omega$ following a procedure similar to \cite{gerhardt:e-r}. 

In a subsequent paper, we discuss boundary regularity. 
If one has boundary regularity, the methods of \cite{athanassenas-finn} can be used to derive higher regularity in smooth domains.

\begin{theorem} \label{Lipschitz}
Let $v$ be a minimiser of  $\J_1$ in $\BV^+(\Omega)$. 
Then $v$ is locally Lipschitz in $\Omega$. 
\end{theorem} 
 
\begin{proof} 
We mollify $v$ over the whole of $\Omega$.  The mollification $v_\eps$ is in
$ C^\infty(\Omega)$, and shares the height bounds derived for $v$ in Section \ref{section4} (that is, bounded above and bounded from below away from zero).  Furthermore, since $v\in\BV(\Omega)$, 
\begin{equation*} 
v_\eps\rightarrow v \text{ in }L^1(\Omega) \text{ and } 
\int_\Omega |Dv_\eps| \rightarrow \int_\Omega |Dv|. 
\end{equation*} 
Corollary \ref{corollary 2.3} for the surface energy and standard convergence results for the remaining energy terms then imply that
\begin{equation}\label{convergence of J1}
\J_1(v_\eps)\rightarrow \J_1(v).
\end{equation}

Let $B\subset\Omega$ be any ball of sufficiently small radius $\rho$, and consider the following two related Dirichlet problems:  
\begin{equation} \label{wrong equation} 
	\begin{cases}\divergenz {\frac{D{w_\eps}}{\sqrt{{w_\eps}^2+|D{w_\eps}|^2}}}=  
{\frac{{w_\eps}}{\sqrt{{w_\eps}^2+|D{w_\eps}|^2}}} +  \ln {w_\eps} \text{ in }B  ,\\  w_\eps=v_\eps \text{ on }\partial B; 
\end{cases}\end{equation} 
and
\begin{equation}\label{wrong equation transformed} 
\begin{cases}\divergenz {\frac{D{u_\eps}}{\sqrt{1+|D{u_\eps}|^2}}}=  
\sigma\left(-{\frac{1}{\sqrt{1+|D{u_\eps}|^2}}} + {u_\eps}  \right)
\text{ in }B,\\  
u_\eps=-\sigma \ln 
(v_\eps) \text{ on }\partial B. 
\end{cases}\end{equation} 
The second expression is in fact a family of problems, indexed by $\sigma\in[0,1]$.   This family is of \emph{mean curvature type}. 
Note that for smooth $w_\eps$ and $u_\eps$, \eqref{wrong equation} is equivalent to \eqref{wrong equation transformed} for $\sigma=1$, with the correspondence $w_\eps=e^{-u_\eps}$.  

Our next step is to solve \eqref{wrong equation transformed} for $\sigma=1$ using the continuity method.  We apply \cite[Theorem 13.8]{GT}.  A prerequisite for this  is to show that a smooth solution 
$u^\sigma$ of \eqref{wrong equation transformed}, for any $\sigma\in[0,1]$, has height and gradient bounds independent of $\sigma$.

The height bound may be found in \cite{Serrin};
 however, the geometric nature of our problem admits a shorter proof which we present as the following lemma.

\begin{lemma} \label{Uniform height estimates}
Let $u^\sigma$ be a  smooth solution to \eqref{wrong equation transformed} corresponding to a $\sigma\in[0,1]$.  Then
\begin{equation}\label{height bound for cont method}
\sup_{B}|u^\sigma|<M_1, 
\end{equation}
where $M_1$ depends only on $\sup_\Omega \left|\ln v_\eps\right|$.
\end{lemma} 
\begin{proof}
We suppose that $u^\sigma$ achieves a positive interior maximum, $u^\sigma(\tilde{x}) = \tilde{M}$ at
some point $\tilde{x} \in B$. If $\tilde{M} > 1$, then the
 mean curvature $H(u^\sigma) = \divergenz {\frac{D{u^\sigma}}{\sqrt{1+
|D{u^\sigma}|^2}}}$ at  $\tilde {x}$ must be strictly positive. 
 But a point of positive mean curvature cannot correspond to
an interior maximum, contradicting the assumption $\tilde{M} > 1$.
We conclude that 
\begin{equation*}
u^\sigma  \le \min \{1,\sup_{\partial B} |\sigma \ln v_\eps| \} \le \min \{1,\sup_{\partial B} |\ln v_\eps| \}.
\end{equation*}
A similar argument shows that $u^\sigma$ has no negative internal minimum,  so $u^\sigma\ge -\sup_{\partial B} |\ln v_\eps|$.
\end{proof}
 
\emph{Continuing the proof of Theorem \ref{Lipschitz}:} 
We find that the gradient bound 
\begin{equation} \label{grad bound for cont method}
\sup_{B}|Du^\sigma|\le M_2
\end{equation} is an application of standard results.  Firstly, an \emph{interior} gradient bound can either be derived by applying a maximum principle to the elliptic equation satisfied by the gradient; or by using \cite[Theorem 4]{ladyzhenskaya-uraltseva}, which gives 
\begin{equation} \label{interior grad bound for usigma}
\sup_{B'}|Du^\sigma|\le M_3
\end{equation} where $B'\subset\subset B$ and $M_3$ is dependent on $\dist(B',\partial B)$, $n$ and $\sup|u^\sigma|$. 

Secondly, a \emph{boundary} gradient estimate 
\begin{equation*}
\sup_{\partial B}|Du^\sigma|\le M_4
\end{equation*} 
results from \cite[Corollary 14.5]{GT} with the structure condition (14.33). Here $M_4$ is dependent on $|\ln v_\eps|_{C^2(\partial B)}$, $n$, $\sup|u^\sigma|$, and $\rho$.    Together these two gradient estimates give us \eqref{grad bound for cont method}.

The conditions for the continuity method being satisfied, the problem \eqref{wrong equation transformed}, with $\sigma=1$, has a $C^{2,\alpha}(B)$ solution which we call $u_\eps$. It has height and gradient bounds \eqref{height bound for cont method} and \eqref{grad bound for cont method}.  
It is also unique: the proof is similar to that of Theorem 2.2 in \cite{athanassenas-finn}, adjusted to Dirichlet boundary data.  

We set $w_\eps=e^{-u_\eps}$.  This is a $C^{2,\alpha}(B)$ solution of \eqref{wrong equation} with height bound $e^{-M_1}\le w_\eps\le e^{M_1}$ and gradient bound $|Dw_\eps|\le M_2e^{M_1}$.    
  
Note that \eqref{wrong equation} 
 is the Euler-Lagrange equation for the energy 
\begin{equation*}
\J_2(w):=\int_B \sqrt{w^2+\left|Dw\right|^2}+ \int_B \int_1^w\
 \ln h \dH \dx, 
\end{equation*}
and so $w_\eps$ is a critical point 
of $\J_2$ in the class of
$H^{1,2}(B)$
 functions with boundary data $v_\eps$.  Furthermore, as the integrand of $\J_2$ is convex in $(w,Dw)$, $w_\eps$ is also a minimiser in this class (see, for example,
the remark in Section 8.2.3 of \cite{Evans}) and hence in the smaller set $C^{2,\alpha}(B)$.

In particular, if we compare $w_\eps$ to $v_\eps$, we have    
\begin{equation} \label{compare energies}
\int_B \sqrt{{w_\eps}^2+\left|Dw_\eps\right|^2}+ \int_B \int_1^{w_\eps}\
 \ln h \dH \dx \le \int_B \sqrt{{v_\eps}^2+\left|Dv_\eps\right|^2}+ \int_B \int_1^{v_\eps}\
 \ln h \dH \dx. 
\end{equation}

Now let $\tilde{v}_\eps$ be defined by 
\begin{equation*}
\tilde{v}_\eps=\begin{cases} w_\eps & \text{ in $B$}\\
v_\eps & \text{ in $\Omega\setminus B$}.
\end{cases}\end{equation*}
Using \eqref{compare energies} for the region $B$ where $\tilde{v}_\eps$ may be different to $v_\eps$, we see that 
%\begin{equation*}
$\J_1(\tilde v_\eps)\le \J_1(v_\eps).$
%\end{equation*}

Now we will show that $\tilde{v}_\eps$ converges to a $\BV(\Omega)$ function which is locally Lipschitz.

Uniform $L^1(\Omega\setminus B)$ bounds are given by the height bounds for $v$ in Section \ref{section4}.  Uniform $L^1(B)$ bounds are given by $\sup_B|w_\eps|\le e^{M_1}$ where $M_1$ is the constant in \eqref{height bound for cont method}; $M_1$ also depends on the height bounds for $v$.

As a consequence of \eqref{convergence of J1}, we may assume that $\J_1(v_\eps)\le \J_1(v)+1$.  Then 
$\J_1(\tilde{v}_\eps)\le \J_1(v)+1$, and so  
\begin{align*}
\int_\Omega \sqrt{\tilde{v}_\eps^2+\left|D\tilde{v}_\eps\right|^2}
&\le \J_1(v)+1-\int_\Omega\int_1^{\tilde{v}_\eps}\ln h \dH\dx+\int_\Sigma 
\beta[1-\tilde{v}_\eps]\ds \\
&\le \J_1(v)+1 +|\Omega| \sup_{\inf v \le h \le \sup v}(h\ln h-h+1)+|\Sigma|\sup_\Sigma |\beta| (1+\sup_\Omega|v|)
\end{align*}
which is bounded above, independently of $\eps$.    Uniform $\BV$ bounds follow as in Lemma \ref{lemma 2.2}.  Therefore a subsequence of $\tilde{v}_\eps$ converges to $v_0\in\BV(\Omega)$, 
and $v_0$ is Lipschitz in $B'$ with bounds given by  \eqref{interior grad bound for usigma}.

Lower semicontinuity of the functional now gives
\begin{equation*}
\J_1(v_0)\le \liminf \J_1(\tilde{v}_\eps)\le \liminf \J_1(v_\eps)= \J_1(v)
\end{equation*}
but as $v$ was assumed to minimise $\J_1$ these must all be equal.  We conclude that there exists a minimiser of $\J_1$ that is locally Lipschitz on interior sets.  
\end{proof} 

Reconsidering the problem of a capillary surface $\S=\graph u$ that minimises the original energy functional $\J$ given in the introduction, we conclude that the found $v$ corresponds to a minimiser in the class $\lbrace w : e^{-w} \in \BV^+(\Omega)\rbrace$.  This solution is given by $u=-\ln v$, and is locally Lipschitz on interior sets.

\def\cprime{$'$}

%\bibliographystyle{acm}

%\bibliography{maria-julie-capillary-bibliography}  
  
\end{document}